\documentclass{amsart}

\usepackage{amsmath, amssymb, amsthm}

\newtheorem{theorem}{Theorem}[section]
\newtheorem{definition}{Definition}[section]

\numberwithin{equation}{section}

\title[Maximal stretching for QC maps]{Existence of quasiconformal maps with maximal stretching on any given countable set}

\author[Tyler Bongers]{Tyler Bongers}
\address{Tyler Bongers, Department of Mathematics, Harvard University, Cambridge, MA 02140}
\email{bongers@math.harvard.edu}

\author[James T. Gill]{James T. Gill}
\address{James T. Gill, Department of Mathematics and Statistics, Saint Louis University, St.\ Louis, MO 63103}
\email{jim.gill@slu.edu}

\date{\today}

\makeatletter
\@namedef{subjclassname@2020}{%
  \textup{2020} Mathematics Subject Classification}
\makeatother

\begin{document}
%\keywords{}

\date{\today}

\begin{abstract}
Quasiconformal maps are homeomorphisms with useful local distortion inequalities; infinitesimally, they map balls to ellipsoids with bounded eccentricity. This leads to a number of useful regularity properties, including quantitative H\"older continuity estimates; on the other hand, one can use the radial stretches to characterize the extremizers for H\"older continuity. In this work, given any bounded countable set in $\mathbb{R}^d$, we will construct an example of a $K$-quasiconformal map which exhibits the maximum stretching at each point of the set. This will provide an example of a quasiconformal map that exhibits the worst-case regularity on a surprisingly large set, and generalizes constructions from the planar setting into $\mathbb{R}^d$.
\end{abstract}

\maketitle

\section{Introduction}
A \textbf{$K$-quasiconformal map} on $\mathbb{R}^d$ is a homeomorphism $F$ whose distributional derivatives satisfy the constraint
\begin{equation}\label{equation:distortion}
\|DF(\vec x)\|^d \le K \det DF(\vec x)
\end{equation}
almost everywhere, where $\|\cdot\|$ stands for the operator norm of the differential matrix $DF$; in order to make sense of the derivatives, we will assume that $F$ lies in the Sobolev class $W^{1, d}(\mathbb{R}^d)$. Morally speaking, quasiconformal maps have nicely bounded distortion properties: the distortion at a point in any one direction (measured by the operator norm) is controlled by the average over all directions (measured by the Jacobian determinant via volume distortion). Infinitesimally, the key geometry is that quasiconformal maps take infinitesimal balls to infinitesimal ellipsoids with bounded eccentricity. These properties arise naturally in elasticity theory and fluid dynamics, where quasiconformal maps have many useful applications; for the general theory, see the classical work of Ahlfors \cite{Ahl06} and the monographs \cite{AstIwaMar09} and \cite{IwaMar01} (handling $\mathbb{C}$ and $\mathbb{R}^n$, respectively).

Quasiconformal maps are natural generalization of conformal maps with less geometric rigidity; however, their geometric properties frequently lead to important regularity results. In the plane, it is a classical result that $K$-quasiconformal maps are $1/K$-H\"older continuous; Mori \cite{Mor56} gave the sharp result
$$|f(z_1) - f(z_2)| \le 16 |z_1 - z_2|^{1/K}$$
for quasiconformal maps $f$ from the unit disk to itself fixing the origin. This generalizes substantially; in fact, it is true in all dimensions that if $F$ is $K$-quasiconformal then
$$|F(\vec y) - F(\vec z)| \le C |\vec y - \vec z|^{1/K}$$
for a constant depending on $F$; see, e.g. \cite[Theorem 3.10.2]{AstIwaMar09} for a proof based on the isoperimetric inequality. (Here and throughout, we will use $| \cdot |$ to stand for the Euclidean norm on $\mathbb{R}^d$). Interpreting this in the context of the distortion estimate \eqref{equation:distortion}, this means that quasiconformal maps can only stretch $\mathbb{R}^d$ so much at any given point.

These results are sharp, in the sense that the exponent cannot be increased: the radial stretch $\vec x \mapsto \vec x |\vec x|^{1/K - 1}$ is $K$-quasiconformal but not $(1/K + \epsilon)$-H\"older continuous for any $\epsilon > 0$. This means that radial stretches exhibit the worst case H\"older continuity among quasiconformal maps; however, at all points except $\vec x = \vec 0$, the map is locally bilipschitz. A natural object of study for a quasiconformal map is therefore the set of points where it exhibits better or worse continuity properties.

There has been a recent push to more precisely quantify how quasiconformal maps stretch or rotate at various points and understand how large the sets can be where a map has bad continuity properties. One way to do this is to define the stretching and rotation exponents of a quasiconformal map in the plane: a quasiconformal map $f : \mathbb{C} \to \mathbb{C}$ is said to stretch with exponent $\alpha$ and rotate with exponent $\gamma$ at $z \in \mathbb{C}$ if there exist scales $r_n \to 0^+$ for which
$$\lim_{n \to \infty} \frac{\log |f(z + r_n) - f(z)|}{\log r_n} = \alpha \text{ and } \lim_{n \to \infty} \frac{\operatorname{arg}(f(z + r_n) - f(z))}{\log|f(z + r_n) - f(z)|} = \gamma,$$
where the argument is interpreted in terms of the winding number of the image of the curve $f([z + r_n, \infty)).$ 

This definition was introduced by Astala, Iwaniec, Prause, and Saksman \cite{AstIwaPraSak15}. The authors classified the exponents attainable at a point by a quasiconformal map and developed upper bounds on the size of the set of points with given stretching and rotation behavior; their bounds were sharp at the level of Hausdorff dimension. A number of authors have sought after sharpness examples, to find maps (whether quasiconformal, or homeomorphisms of finite distortion, or other useful classes) which exhibit the worst-case stretching and rotation; see e.g. the work of Clop, Hitruhin, and Sengupta \cite{CloHitSen21} for the worst-case rotation of a homeomorphism with $L^p$ distortion function.

Still working in the plane, the first named author showed in \cite[Theorem 3.4]{Bon19} additional sharpness examples for stretching in the extreme case $\alpha = 1/K$: given any countable set $\Lambda$ contained in the unit disk, it is possible to construct a $K$-quasiconformal map which exhibits the worst-case stretching exponent at every point in $\Lambda$.

The goal of this paper is to extend this result to higher dimensional quasiconformal maps; that is, to show that given a bounded countable set $\Lambda$, there is a $K$-quasiconformal map that exhibits the worst-case stretching behavior at all points of $\Lambda$. To this end, we need to precisely define the stretching exponent:
\begin{definition}
A map $F : \mathbb{R}^d \to \mathbb{R}^d$ is said to \emph{stretch with exponent $\alpha$} at the point $\vec x$ if there is a sequence of scales $r_n \to 0^+$ and a sequence of unit vectors $\vec u_n$ for which
$$\lim_{n \to \infty} \frac{\log |F(\vec x + r_n \vec u_n) - F(\vec x)|}{\log r_n} = \alpha.$$
If such a sequence exists, we will call $\alpha$ a \emph{stretching exponent} of $F$ at $\vec x$.
\end{definition}

With this notation, we can now state the main result of the paper:
\begin{theorem}\label{theorem:existence}
Fix a bounded countable set $\Lambda = \{\lambda_1, \lambda_2, ...\} \subseteq \mathbb{R}^d$ and a $K > 1$. There is a $K$-quasiconformal map $F_{\Lambda}$ which stretches with exponent $1/K$ at every point in $\Lambda$.
\end{theorem}

This is a direct generalization of \cite[Theorem 3.4]{Bon19}, although the proof here is substantially more complicated. In the complex setting, one considers a sum of radial stretches, and builds a map similar to
$$f(z) = \sum_{n = 1}^{\infty} \frac 1 {2^n} (z - \lambda_n)|z - \lambda_n|^{1/K}.$$
It is quite rare that the sum of two quasiconformal maps is quasiconformal, let alone injective. However, there is a key positivity property that can be used to show that $f$ is, in fact, quasiconformal: it is true that
\begin{equation}\label{equation:complex_positivity}
\frac{\partial}{\partial z} z |z|^{1/K - 1} \ge 0
\end{equation}
for all $z \ne 0$. It turns out that this is exactly what is needed to show that $f$ satisfies a Beltrami equation, which in turn can be bootstrapped into showing it is quasiconformal.

In the higher dimensional, real setting there are a number of obstacles which have to be overcome to make an analogous construction. First, one must identify what notion of ``positivity'' is relevant here and which can replace \eqref{equation:complex_positivity}; in this case, it will be the fact that the differential matrix of a radial stretch is positive definite. Secondly, and more crucially, the distortion inequality \eqref{equation:distortion} involves the operator norm $\|DF\|$ and the Jacobian $\det DF$, both of which are very nonlinear objects. This requires substantial analysis of the Jacobian of the sum of radial stretches, involving detailed work with the spectra of various matrices. After these complications are addressed, we are then able to show that our map is quasiconformal and exhibits the correct stretching at all points of $\Lambda$.

\section{Proof of Theorem \ref{theorem:existence}}
%We now turn to the proof of Theorem \ref{theorem:existence}, breaking the argument into several stages.
\subsection{Construction of $F_{\Lambda}$} The basic building block of $F_{\Lambda}$ will be the radial stretching map located at the origin:
$$S(\vec x) = \vec x |\vec x|^{1/K - 1}.$$
This is a $K$-quasiconformal map from $\mathbb{R}^d$ to $\mathbb{R}^d$ and stretches with exponent $1/K$ at the origin (and is bilipschitz on any compact set separated from the origin). We will need a detailed study of its differential matrix; fix a point $\vec x$ within the unit ball. The $(i,j)$ component of the matrix $DS$ is given by
\begin{align*}
[DS]_{i,j} &= \frac{\partial}{\partial x_j} (x_1^2 + \cdots + x_d^2)^{\frac{1 - K}{2K}} x_i \\
&= \frac{1 - K}{2K} (x_1^2 + \cdots + x_d^2)^{\frac{1 - K}{2K} - 1} \cdot 2x_j \cdot x_i + (x_1^2 + \cdots x_d^2)^{\frac{1 - K}{2K}} \frac{\partial x_i}{\partial x_j} \\
&= \frac{1 - K}{K} |\vec x|^{\frac{1 - K}{K} - 2} x_i x_j + |\vec x|^{\frac{1 - K}{K}} \delta_{i, j} \\
&= |\vec x|^{1/K - 1} \left(\delta_{i, j} - \left(1 - \frac 1 K\right) \frac{x_i x_j}{|\vec x|^2}\right),
\end{align*}
where $\delta_{i, j} = 1$ if $i = j$ and is zero otherwise. This gives a particularly special form to the matrix $DS$: it is essentially a perturbation of a (properly scaled) identity matrix via a rank one operator. To be specific, we have
\begin{equation}\label{equation:stretch_matrix}
DS(\vec x) = |\vec x|^{1/K - 1} \left(I - \left(1 - \frac 1 K\right) \frac{\vec x \vec x^T}{|\vec x|^2}\right).
\end{equation}

We need to establish a few properties of this matrix before continuing; note that $DS$ is a symmetric matrix. Next, and most important, is that this matrix is positive definite. To see this, observe that if $\vec w \in \mathbb{R}^d$, we have that
\begin{align*}
\vec w^T (\vec x \vec x^T) \vec w &= (\vec x^T \vec w)^T (\vec x^T \vec w) = |\vec x^T \vec w|^2
\end{align*}
is a nonnegative real number at most $|\vec x|^2 \cdot |\vec w|^2$ via the Cauchy-Schwarz inequaltiy. This implies that
$$\left|\vec w^T \left(\left(1 - \frac 1 K\right) \frac{\vec x \vec x^T}{|\vec x|^2} \right) \vec w\right| \le \left(1 - \frac 1 K\right) |\vec w|^2 < |\vec w|^2$$
since $K > 1$. The positive definiteness of $DS$ on the unit ball follows immediately from this upper bound. 

We are now ready to construct the map $F_{\Lambda}$ from slightly modified versions of $S$. For a fixed point $\lambda$, define
$$S_{\lambda}(\vec x) = (\vec x - \lambda) |\vec x - \lambda|^{1/K - 1}$$
and set
\begin{equation}\label{equation:function_definition}
F_{\Lambda}(\vec x) = \sum_{n = 1}^{\infty} \frac 1 {2^n} S_{\lambda_n}(\vec x).
\end{equation}
Note that by translating \eqref{equation:stretch_matrix}, 
$$DS_{\lambda}(\vec x) = |\vec x - \lambda|^{1/K -1 } \left(I - \left(1 - \frac 1 K\right) \frac{(\vec x - \lambda)(\vec x - \lambda)^T}{|\vec x - \lambda|^2}\right).$$
Furthermore, because the series rapidly converges, $F_{\Lambda}$ lies in the Sobolev class $W^{1, d}(\mathbb{R}^d)$. One way to argue this is to observe that $F_{\Lambda}$ is absolutely continuous along every vertical and horizontal line which does not intersect $\Lambda$ and study the classical difference quotients used to define its partial derivatives. In any case, we find that
$$DF_{\Lambda}(x) = \sum_{n = 1}^{\infty} \frac 1 {2^n} DS_{\lambda_n}(\vec x).$$
Each $DS_{\lambda_n}$ is symmetric and positive definite, being just a translate of $DS$; therefore, $DF_{\Lambda_n}$ is also symmetric and positive definite and this shows that $F_{\Lambda}$ is injective. There are now two aspects left to verify: that $F_{\Lambda}$ satisfies the distortion inequality for $K$-quasiconformal maps and that it exhibits the correct stretching exponent on $\Lambda$. 

\subsection{Establishing $K$-quasiconformality} 
In order to establish the distortion inequality \eqref{equation:distortion}, we rewrite the differential matrix with some extra notation. Starting with
$$DF_{\Lambda}(\vec x) = \sum_{n = 1}^{\infty} \frac 1 {2^n} |\vec x - \lambda_n|^{1/K - 1} \left(I - \left(1 - \frac 1 K\right) \frac{(\vec x - \lambda_n)(\vec x - \lambda_n)^T}{|\vec x - \lambda_n|^2}\right),$$
set $\vec w_n = (\vec x - \lambda_n) / |\vec x - \lambda_n|$, $W(\vec x) = \sum_{n = 1}^{\infty} 2^{-n} |\vec x - \lambda_n|^{1/K - 1}$, and
$$\eta_n(\vec x) = \frac{2^{-n} |\vec x - \lambda_n|^{1/K - 1}}{W(x)},$$
defined for all $\vec x \notin \Lambda$. Note that for such $\vec x$, we have $\sum_{n = 1}^{\infty} \eta_n(\vec x) = 1$. We then have
\begin{align*}
DF_{\Lambda}(\vec x) &= \left(\sum_{n = 1}^{\infty} \frac 1 {2^n} |\vec x - \lambda_n|^{1/K - 1} \right) I - \left(1 - \frac 1 K\right) \sum_{n = 1}^{\infty} \frac 1 {2^n} |\vec x - \lambda_n|^{1/K -1 } \vec w_n \vec w_n^T \\
&= W(\vec x) I - \left(1 - \frac 1 K\right) W(\vec x) \sum_{n = 1}^{\infty} \frac{2^{-n} |\vec x - \lambda_n|^{1/K -1 }}{W(\vec x)} \vec w_n \vec w_n^T \\
&= W(\vec x) I - \left(1 - \frac 1 K\right) W(\vec x) \sum_{n = 1}^{\infty} \eta_n(\vec x) \vec w_n \vec w_n^T \\
&= W(\vec x) \left(I - \left(1 - \frac 1 K\right) \sum_{n = 1}^{\infty} \eta_n(\vec x) \vec w_n \vec w_n^T\right).
\end{align*}
For a final set of notation, set $\alpha = 1 - 1/K \in (0, 1)$ and
$$B(\vec x) = \sum_{n = 1}^{\infty} \eta_n(\vec x) \vec w_n \vec w_n^T,$$
noting that $B(\vec x)$ is a matrix defined for each $\vec x \notin \Lambda$; in this new notation, $DF_{\Lambda}(\vec x) = W(\vec x) (I - \alpha B(\vec x)).$

Our aim is to show that $\|DF_{\Lambda}\|^{d} \le K \det(DF_{\Lambda});$ since both sides are homogeneous of the same degree and $W(\vec x)$ is a scalar function, it is sufficient to establish that for each $\vec x \notin \Lambda$ we have 
\begin{equation}\label{equation:distortion_inequality}
\left\|I - \alpha B(\vec x)\right\|^{d} \le K \det\left(I - \alpha B(\vec x)\right).
\end{equation}
Since $\Lambda$ has Lebesgue measure zero, this is sufficient to prove quasiconformality. To this end, we will show that the operator norm is at most $1$ and the determinant is at least $1/K$.

In preparation for this, we need to analyze the spectrum of the matrix $I - \alpha B(\vec x)$, starting with the symmetric matrix $B(\vec x)$. If $y \in \mathbb{R}^d$ is a unit vector,
\begin{align*}
\vec y^T B(\vec x) \vec y &= \sum_{n = 1}^{\infty} \eta_n(\vec x) \vec y^T \vec w_n \vec w_n^T \vec y \\
&= \sum_{n = 1}^{\infty} \eta_n(\vec x) |\vec w_n^T \vec y|^2. 
\end{align*}
Estimating $|\vec w_n^T \vec y| \le |\vec w_n| \cdot |\vec y| \le 1$ and recalling that $\sum_n \eta_n(\vec x) = 1$, we have
$$\vec y^T B(\vec x) \vec y \le 1$$
for all unit vectors $\vec y$. On the other hand, each $\eta_n(\vec x)$ is nonnegative and so $\vec y^T B(\vec x) \vec y \ge 0$. Therefore, $B(\vec x)$ is a positive semidefinite matrix with operator norm at most $1$. This implies that the spectrum $\sigma(B(\vec x))$ consists of $d$ eigenvalues $\{\sigma_1, ..., \sigma_d\} \subseteq [0, 1]$. 

There is another important fact about the spectrum of $B(\vec x)$ which we will need: the sum of the eigenvalues is $1$. To see this, note that the trace of $B(\vec x)$ is
$$\operatorname{Tr} B(\vec x) = \sum_{n = 1}^{\infty} \eta_n(\vec x) \operatorname{Tr} (\vec w_n \vec w_n^T) = \sum_{n = 1}^{\infty} \eta_n(\vec x) |\vec w_n|^2 = 1.$$
after using that each $w_n$ is a unit vector. Recalling that the trace of a matrix is the sum of its eigenvalues, the claim follows.

We also need to relate the spectrum of $B(\vec x)$ to that of $I - \alpha B(\vec x)$. When $\vec v_i$ is an eigenvector of $B(\vec x)$ associated to the eigenvalue $\sigma_i$, it follows that $\vec v_i$ is also an eigenvector of $I - \alpha B(\vec x)$ but now with eigenvalue $1 - \alpha \sigma_i$. It follows that the spectrum of $I - \alpha B(\vec x)$ is the set
$$\sigma(I - \alpha B(\vec x)) = \{1 - \alpha \sigma_1, \dots, 1 - \alpha \sigma_d\}.$$
Because each $\sigma_i \in [0, 1]$ and $\alpha \in (0, 1)$, it follows that the eigenvalues of $I - \alpha B(\vec x)$ all lie in the interval $(0, 1]$. 

We are now ready to bound the relevant quantities for the distortion inequality \eqref{equation:distortion_inequality}. Since $I - \alpha B(\vec x)$ is a real symmetrix matrix, it is orthogonally diagonalizable; all its eigenvalues are at most $1$, and so
\begin{equation}\label{equation:operator_norm}
\|I - \alpha B(\vec x)\| \le \max \{\sigma_1, ..., \sigma_d\} \le 1.
\end{equation}
It therefore remains to bound the determinant from below; to this end, write the determinant as the product of the eigenvalues:
\begin{align}\label{equation:determinant_bound}
\det(I - \alpha B(\vec x)) &= \prod_{i = 1}^d (1 - \alpha \sigma_i) \nonumber \\
&= \sum_{k = 0}^d (-\alpha)^k \sum_{i_1 < ... < i_k} \sigma_{i_i} \cdots \sigma_{i_k} \nonumber \\
&= 1 - \alpha (\sigma_1 + \cdots + \sigma_d) + \sum_{k = 2}^d (-\alpha)^k \sum_{i_1 < ... < i_k} \sigma_{i_1} \cdots \sigma_{i_k} \nonumber \\
&= \frac 1 K + \sum_{k = 2}^d (-\alpha)^k \sum_{i_1 < ... < i_k} \sigma_{i_1} \cdots \sigma_{i_k},
\end{align}
where we have used the fact that the trace of $B(\vec x)$ is $1$ together with the definition of $\alpha$ in the final line. We must show that the summation is nonnegative; the obstacle is its alternating nature. For notation, define the symmetric polynomials 
$$P_k(\sigma_1, ..., \sigma_d) = \sum_{i_1 < ... < i_k} \sigma_{i_1} \cdots \sigma_{i_k}.$$
Recalling that $\sigma_1 + \cdots + \sigma_d = 1$ and that each $\sigma_i \ge 0$, we have
\begin{align*}
P_k(\sigma_1, ..., \sigma_d) &= P_k(\sigma_1, ..., \sigma_d)\cdot(\sigma_1 + \cdots + \sigma_d) \\
&= \sum_{j = 1}^d \sum_{i_1 < ... < i_k} \sigma_{i_1} \cdots \sigma_{i_k} \sigma_j \\
&= (k + 1) \sum_{j_1 < ... < j_{k + 1}} \sigma_{j_1} \cdots \sigma_{j_{k + 1}} + \sum_{\text{duplicated indices}} \sigma_{\ell_1} \cdots \sigma_{\ell_{k + 1}} \\
&\ge (k + 1) P_{k + 1}(\sigma_1, ..., \sigma_d).
\end{align*}
The third line follows from the fact that each $(k + 1)$-tuple can be decomposed as a $k$-tuple and a single index in $(k + 1)$ ways, and the fourth line uses the definition of $P_{k + 1}$ together with the fact that every summand is nonnegative. In fact, we only need the weaker inequality $P_k(\sigma_1, ..., \sigma_d) \ge P_{k + 1}(\sigma_1, ..., \sigma_d)$; since $\alpha \in (0, 1)$ and the first term is nonnegative, it follows that
$$\sum_{k = 2}^d (-\alpha)^k P_k(\sigma_1, ..., \sigma_d) \ge 0.$$
Inserting this into \eqref{equation:determinant_bound} shows that $\det(I - \alpha B(\vec x)) \ge \frac 1 K$, whose combination with with \eqref{equation:operator_norm} gives
$$\|I - \alpha B(\vec x)\|^{d} \le 1 \le K \det(I - \alpha B(\vec x))$$
for all $\vec x \notin \Lambda$. This yields the distortion inequality \eqref{equation:distortion_inequality} and completes the verification that $F_{\Lambda}$ is $K$-quasiconformal.

\subsection{Stretching exponent for points in $\Lambda$} All that remains now is to study the stretching behavior on $\Lambda$; this will follow the lines of the argument in \cite{Bon19}, but there are complications in working on $\mathbb{R}^d$ rather than $\mathbb{C}$. To begin, fix an index $N$; we will show that for all sufficiently small $r$ that
\begin{equation}\label{equation:stretching_goal}
|F_{\Lambda}(\lambda_N + r e_1) - F_{\Lambda}(\lambda_N)| \gtrsim r^{1/K}.
\end{equation}
After translations, we may assume that $\lambda_N = \vec 0$. Returning to the definition \eqref{equation:function_definition}, we have
\begin{align*}
F_{\Lambda}(r \vec e_1) - F_{\Lambda}(\vec 0) &= \sum_{n = 1}^{\infty} \frac{1}{2^n} \left(S_{\lambda_n}(r \vec e_1) - S_{\lambda_n}(\vec 0)\right) \\
&= \frac{1}{2^N} (r\vec e_1) |r \vec e_1|^{1/K - 1} - \sum_{n \ne N} \frac{1}{2^n} \left(S_{\lambda_n}(r\vec e_1) - S_{\lambda_n}(\vec 0)\right) \\
&= \frac{r^{1/K}}{2^N} \vec e_1 - \sum_{n \ne N} \frac{1}{2^n} \left(S_{\lambda_n}(r\vec e_1) - S_{\lambda_n}(\vec 0)\right).
\end{align*}
We will show that if $r$ is sufficiently close to zero, then the latter term is negligible in comparison to the first term. 

Roughly speaking, there are two kinds of points $\lambda_n$ which we will consider. First, there are points with ``small'' indices $n$; these points can be treated as far away from $\lambda_N = \vec 0$, in which case we use the Lipschitz nature of $S$. Secondly, there are points with ``large'' indices $n$ which may be arbitrarily close to $\lambda_N$; in this case, we use the exponential decay and a basic H\"{o}lder continuity estimate. We will take the cutoff point to an index $N^*$ to be specified later. We therefore need to establish a pair of inequalities: for each $\epsilon > 0$, there exists an $r^* > 0$ such that whenever $0 < r < r^*$ we have
\begin{equation}\label{equation:stretching_goal_far}
\sum_{\substack{1 \le n < N^* \\ n \ne N}} \frac 1 {2^n} \left|S_{\lambda_n}(r\vec e_1) - S_{\lambda_n}(\vec 0)\right| \le \epsilon r^{1/K}
\end{equation}
and
\begin{equation}\label{equation:stretching_goal_near}
\sum_{n= N^*}^{\infty} \frac{1}{2^n} \left|S_{\lambda_n}(r\vec e_1) - S_{\lambda_n}(\vec 0)\right| \le \epsilon r^{1/K}.
\end{equation}

In order to establish these inequalities, we need to study the precise behavior of the stretching maps $S_{\lambda}$. We begin with studying the stretch in a particular direction adapted to $\lambda$: 
\begin{align}\label{equation:radial_change}
\left|S_{\lambda} \left(r \frac{\lambda}{|\lambda|}\right) - S_{\lambda}(\vec 0)\right| &= \left|\left(r \frac{\lambda}{|\lambda|} - \lambda\right) \left|r \frac{\lambda}{|\lambda|} - \lambda\right|^{1/K - 1} - (-\lambda) |-\lambda|^{1/K - 1}\right| \nonumber \\
&= \left|\frac{\lambda}{|\lambda|} (r - |\lambda|) \cdot \left|\frac{\lambda}{|\lambda|} (r - |\lambda|)\right|^{1/K - 1} + \frac{\lambda}{|\lambda|} |\lambda|^{1/K} \right| \nonumber \\
&= \left|(r - |\lambda|) \cdot |r - |\lambda||^{1/K - 1} + |\lambda|^{1/K} \right| \nonumber \\
&= \left||\lambda| \left(\frac{r}{|\lambda|} - 1\right) \cdot |\lambda|^{1/K - 1} \left|\frac{r}{|\lambda|} - 1\right|^{1/K - 1} + |\lambda|^{1/K}\right| \nonumber \\
&= |\lambda|^{1/K} \left|\left(\frac{r}{|\lambda|} - 1\right)\left|\frac{r}{|\lambda|} - 1\right|^{1/K - 1} + 1\right| \nonumber \\
&= |\lambda|^{1/K} \left|1 + \frac{r/|\lambda| - 1}{|r/|\lambda| - 1|} \cdot |r/|\lambda| - 1|^{1/K}\right|.
\end{align}
Set $t = r/|\lambda| \in (0, \infty)$, in which case we have
$$f(t) := \left|S_{\lambda} \left(r \frac{\lambda}{|\lambda|}\right) - S_{\lambda}(0)\right| = |\lambda|^{1/K} \left|1 + \frac{t - 1}{|t - 1|}|t - 1|^{1/K}\right|.$$
We now separate into two cases based on the relative scale of $r$ and $|\lambda|$:
\begin{itemize}
\item If $t \ge 1$, then 
$$f(t) = |\lambda|^{1/K} \left(|t - 1|^{1/K} + 1\right) \lesssim |\lambda|^{1/K} t^{1/K}.$$
\item If $t \in (0, 1)$, then
$$f(t) = |\lambda|^{1/K} (1 - (1 - t)^{1/K}).$$
Taking the Taylor expansion around $0$, we have
$$1 - (1 - t)^{1/K} = 1 - \left[1 - \frac 1 K t + O(t^2)\right] = \frac 1 K t + O(t^2).$$
This leads to $f(t) \lesssim |\lambda|^{1/K} t$. 
\end{itemize}
In each case, we have that $f(t) \lesssim \min\{t, t^{1/K}\}.$ Returning to \eqref{equation:radial_change}, we have established that
\begin{equation}\label{equation:stretch_special_direction}
\left|S_{\lambda} \left(r \frac{\lambda}{|\lambda|}\right) - S_{\lambda}(\vec 0)\right| \lesssim |\lambda|^{1/K} \min\left\{\frac{r}{|\lambda|}, \left(\frac{r}{|\lambda|}\right)^{1/K}\right\}.
\end{equation}
In order to establish the inequalities \eqref{equation:stretching_goal_far} and \eqref{equation:stretching_goal_near}, we need a slight extension of this, in that the unit vector $\lambda / |\lambda|$ needs to be replaced by an arbitrary unit vector $\vec u$. That is, we need to show that
\begin{equation}\label{equation:stretch_all_directions}
|S_{\lambda}(r \vec u) - S_{\lambda}(\vec 0)| \lesssim |\lambda|^{1/K} \min \left\{\frac{r}{|\lambda|}, \left(\frac{r}{|\lambda|}\right)^{1/K}\right\}.
\end{equation}
However, this follows immediately from the fact that a (global) quasiconformal map on $\mathbb{R}^d$ is also quasisymmetric.

We may now complete the analysis by selecting the radius scale $r^*$ and the cutoff $N^*$. Recall that by \eqref{equation:stretch_all_directions}, there is a constant $C < \infty$ independent of $\lambda$ for which
\begin{equation}\label{equation:alternatives}
\left|S_{\lambda} (r \vec e_1) - S_{\lambda}(\vec 0)\right| \le C \min\{r^{1/K}, r \cdot |\lambda|^{1/K - 1}\}.
\end{equation}
Set the cutoff $N^* = N + A$, where $A$ is large enough that $C / 2^{N + A - 1} < \epsilon$; in this case, we use the first alternative of \eqref{equation:alternatives} to get
\begin{align*}
\sum_{n = N^*}^{\infty} \frac{1}{2^n} |S_{\lambda_n}(r \vec e_1) - S_{\lambda_n}(\vec 0)| &\le \sum_{n = N + A}^{\infty} \frac{1}{2^n} \cdot C r^{1/K} = \frac{2C}{2^{N + A}} r^{1/K} < \epsilon r^{1/K}.
\end{align*}
This establishes \eqref{equation:stretching_goal_far} as desired.

Finally, we need to handle small indices as well. For this, use the linear estimate available in \eqref{equation:alternatives}:
\begin{align}\label{equation:radius_selection}
\sum_{\substack{1 \le n \le N^* \\ n \ne N}} \frac{1}{2^n} |S_{\lambda_n}(r \vec e_1) - S_{\lambda_n}(\vec 0)| &\le \sum_{\substack{1 \le n \le N^* \\ n \ne N}} \frac{1}{2^n} C r |\lambda_n|^{1/K - 1} \nonumber \\
&\le C r \max_{\substack{1 \le n \le N^* \\ n \ne N}} |\lambda_n|^{1/K - 1} \nonumber \\
&= C r^{1/K} \max_{\substack{1 \le n \le N^* \\ n \ne N}} \left(\frac r {|\lambda_n|}\right)^{1 - 1/K} \nonumber \\
&= C r^{1/K} \left(\frac{r}{\min\limits_{\substack{1 \le n \le N^* \\ n \ne N}} |\lambda_n|}\right)^{1 - 1/K}.
\end{align}
Note that since $n \ne N$ and $N^* < \infty$, there is a uniform, positive lower bound on $|\lambda_n|$ for this context. Let $\rho > 0$ be this minimum radius and set 
$$r^* = \rho \left(\frac{\epsilon}{C}\right)^{1 / (1 - 1/K)}.$$
Note that the selection of $\rho$ depends only on $\epsilon$, the set $\Lambda$, and the value $N^*$; but $N^*$ depends only on the constant in inequality \eqref{equation:stretch_all_directions}, which ultimately only depends on $K$. That is, the selection of $r^*$ depends only on $\epsilon$, $\Lambda$, and $K$ and is independent of $N$. Inequality \eqref{equation:stretching_goal_near} now follows immediately from inserting the selection of $r^*$ into \eqref{equation:radius_selection}, completing the analysis of the stretching at $\Lambda_N$ and proving \eqref{equation:stretching_goal}. Since $N \in \mathbb{N}$ was arbitrary, the proof is completed.

\bibliographystyle{plain}
\bibliography{qc_bib}

\end{document}